\documentstyle[12pt,amsfonts]{article}

\begin{document}
\title{H\"{o}lder continuity of solutions of
supercritical dissipative hydrodynamic transport equations}

\author{Peter Constantin
\\Department of Mathematics \\University of Chicago
\\ 5734 S. University Avenue\\Chicago, IL 60637 \\ \\E-mail: const@cs.uchicago.edu
\\  \\ \\Jiahong Wu
\\Department of Mathematics\\ Oklahoma State University\\ Stillwater, OK 74078\\ \\
E-mail: jiahong@math.okstate.edu}
\date{}
\maketitle

\newtheorem{thm}{Theorem}[section]
\newtheorem{thmA}[thm]{Theorem A}
\newtheorem{cor}[thm]{Corollary}
\newtheorem{prop}[thm]{Proposition}
\newtheorem{define}[thm]{Definition}
\newtheorem{rem}[thm]{Remark}
\newtheorem{example}[thm]{Example}
\newtheorem{lemma}[thm]{Lemma}
\def\theequation{\thesection.\arabic{equation}}

\vspace{.15in} \noindent {\bf Abstract}. We examine the
regularity of weak solutions of quasi-geostrophic (QG) type
equations with supercritical ($\alpha <1/2$) dissipation $(-\Delta)^\alpha$.
This study is motivated by a recent work of
Caffarelli and Vasseur, in which they study the global regularity
issue for the critical ($\alpha = 1/2$) QG equation \cite{CV}.  Their approach
successively increases the regularity levels of Leray-Hopf weak
solutions: from $L^2$ to $L^\infty$, from $L^\infty$ to
H\"{o}lder ($C^{\delta}$, $\delta>0$), and
from H\"{o}lder to classical solutions. In the supercritical case,
Leray-Hopf weak solutions can still be shown to be $L^\infty$, but
it does not appear that their approach can be easily extended to
establish the H\"{o}lder continuity of $L^\infty$ solutions. In
order for their approach to work, we require the velocity to be
in the H\"{o}lder space $C^{1-2\alpha}$. Higher regularity starting
from $C^\delta$ with $\delta>1-2\alpha$ can be established through
Besov space techniques and will be presented elsewhere \cite{CW6}.

\vspace{.2in} \noindent {\bf AMS (MOS) Numbers:} 76D03, 35Q35

\vspace{.1in} \noindent {\bf Keywords:} the dissipative
quasi-geostrophic equation, regularity, supercritical dissipation,
weak solutions.

\newpage
\section{Introduction}
\setcounter{equation}{0} \label{sec:1}

This paper studies the regularity of Leray-Hopf weak solutions of
the dissipative QG equation of the form
\begin{equation}\label{QG}
\left\{
\begin{array}{l}
\partial_t \theta + u\cdot\nabla \theta +\kappa (-\Delta)^\alpha\theta=
0,\quad x\in {\Bbb R}^n,\,\,t>0,\\ \\
u={\cal R}(\theta),\quad\nabla \cdot u=0,\quad x\in {\Bbb
R}^n,\,\,t>0,
\end{array}
\right.
\end{equation}
where $\theta=\theta(x,t)$ is a scalar function, $\kappa>0$ and
$\alpha>0$ are parameters, and ${\cal R}$ is a standard singular
integral operator. The fractional Laplace operator
$(-\Delta)^\alpha$ is defined through the Fourier transform
$$
\widehat{(-\Delta)^\alpha f}(\xi) = |\xi|^{2\alpha}
\widehat{f}(\xi),\quad \xi\in {\Bbb R}^n.
$$
(\ref{QG}) generalizes the 2D dissipative QG equation (see
\cite{Con},\cite{CMT},\cite{HPGS},\cite{Pe} and the references
therein). The main  mathematical question concerning the 2-D
dissipative QG equation is whether or not it has a global in time
smooth solution for any prescribed smooth initial data. In the
subcritical case $\alpha>\frac{1}{2}$, the dissipative QG equation
has been shown to possess a unique global smooth solution for every
sufficiently smooth initial data (see \cite{CW5},\cite{Res}). In
contrast, when $\alpha < \frac12$, the question of global existence
is still open. Recently this problem has attracted a significant
amount of research
(\cite{CV},\cite{Cha},\cite{ChL},\cite{CMZ},\cite{Con},\cite{CCW},\cite{CC},\cite{HK}
\cite{Ju},\cite{KNV},\cite{MarLr},\cite{Sch},\cite{Sch2},\cite{Wu02},\cite{Wu3},
\cite{Wu4},\cite{Wu77}). In Constantin, C\'{o}rdoba and Wu
\cite{CCW}, we proved in the critical case ($\alpha=\frac12$) the
global existence and uniqueness of classical solutions corresponding
to any initial data with $L^\infty$-norm comparable to or less than
the diffusion coefficient $\kappa$. In a recent work \cite{KNV},
Kiselev, Nazarov and Volberg proved that smooth global solutions
persist for any $C^\infty$ periodic initial data \cite{CCW}, for the
critical QG equation. Also recently, Caffarelli and Vasseur
\cite{CV} proved the global regularity of the Leray-Hopf weak
solutions to the critical QG equation in the whole space.

\vspace{.1in} We focus our attention on the supercritical case
$\alpha<\frac{1}{2}$. Our study is motivated by the work of
Caffarelli and Vasseur in the critical case. Roughly speaking, the
Caffarelli-Vasseur approach consists of three main steps. The first
step shows that a Leray-Hopf weak solution emanating from an initial
data $\theta_0\in L^2$ is actually in $L^\infty({\Bbb
R}^n\times(0,\infty))$. The second step proves that the
$L^\infty$-solution is $C^\gamma$-regular, for some $\gamma>0$. For
this purpose, they represent the diffusion operator $\Lambda\equiv
(-\Delta)^{1/2}$ as the normal derivative of the harmonic extension
$L$ from $C_0^\infty({\Bbb R}^n)$ to $C^\infty_0({\Bbb
R}^n\times{\Bbb R}^+)$ and then exploit a version of De Giorgi's
isoperimetric inequality to prove the H\"{o}lder continuity. The
third step improves the H\"{o}lder continuity to $C^{1,\beta}$, the
regularity level of classical solutions.

\vspace{.1in} We examine the approach of Caffarelli and Vasseur to
see if it can be extended to the super-critical case. The first step
of their approach can be modified to suit the supercritical case:
any Leray-Hopf weak solution can still be shown to be $L^\infty$ for
any $x\in {\Bbb R}^n$ and $t>0$ (see Theorem \ref{l8bound}).
Corresponding to their third step, we can show that any weak solution
already in the H\"{o}lder class $C^\delta$ with $\delta>1-2\alpha$, is
actually a global classical solution. This result is
established by representing the H\"{o}lder space functions in terms
of the Littlewood-Paley decomposition and using Besov space
techniques. We will present this result in a separate paper
\cite{CW6}. We do not know if any solution in H\"{o}lder space
$C^\gamma$ with arbitrary $\gamma>0$ is smooth, and therefore there exists
a significant potential obstacle to the program:  even if
all Leray-Hopf solutions are $C^\gamma$, $\gamma>0$, it may still
be the case that only those solutions for which $\gamma>1-2\alpha$
are actually smooth.
If this would be true, then the critical case would be a fortuitous
one, ($1-2\alpha = 0$). If, however, all Leray-Hopf solutions are smooth,
then providing a proof of this fact would require a new idea.

\vspace{.1in} The most challenging part is how to establish the
H\"{o}lder continuity of the $L^\infty$-solutions. It does not
appear that the approach of Caffarelli and Vasseur can be easily
extended to the supercritical case. In the critical case, Caffarelli
and Vasseur lifted $\theta$ from ${\Bbb R}^n$ to a harmonic function
$\theta^*$ in the upper-half space ${\Bbb R}^n\times {\Bbb R}^+$
with boundary data on ${\Bbb R}^n$ being $\theta$. The fractional
derivative $(-\Delta)^{\frac12}\theta$ is then expressed as the
normal derivative of $\theta^*$ on the boundary ${\Bbb R}^n$ and the
$\mathring{H}^1$-norm of $\theta^*$ is then bounded by the natural
energy of $\theta$. Taking the advantage of the nice properties of
harmonic functions, they were able to obtain a diminishing
oscillation result for $\theta^*$ in a box near the origin. More
precisely, if $\theta^*$ satisfying $|\theta^*| \le 2$ in the box,
then $\theta^*$ satisfies in a smaller box centered at the origin
$$
\sup \,\theta^* - \inf\, \theta^* <4-\lambda^*
$$
for some $\lambda^*>0$. The proof of this result relies on a local
energy inequality, an isoperimetric inequality of De Giorgi and two
lengthy technical lemmas. Examining the proof reveals that
$\lambda^*$ depends on the BMO-norm of the velocity $u$. To show the
H\"{o}lder continuity at a point, they zoom in at this point by
considering a sequence of functions $\theta^*_k$ and $u_k$ with
$(\theta_k, u_k)$ satisfying the critical QG equation. This process
is carried out through the natural scaling invariance that
$(\theta(\mu x, \mu t), u(\mu x, \mu t))$ solves the critical QG
equation if $(\theta,u)$ does so. Applying the diminishing
oscillation result to this sequence leads to the H\"{o}lder
continuity of $\theta^*$. An important point is that the BMO-norm of
$u_k$ is preserved in this scaling process.

\vspace{.1in} In the supercritical case, the diminishing oscillation
result can still be established by following the idea of Caffarelli
and Vasseur (see Theorem \ref{dios}). However,  the scaling
invariance is now represented by $\mu^{2\alpha-1}\theta(\mu x,
\mu^{2\alpha} t)$ and $\mu^{2\alpha-1}u(\mu x, \mu^{2\alpha} t)$ and
the BMO-norm deteriorates every time the solution is rescaled. This
is where the approach of Caffarelli and Vasseur stops working for
the supercritical case. If we make the assumption that $u\in
C^{1-2\alpha}$, then the scaling process preserves this norm and we
can still establish the H\"{o}lder continuity of $\theta$. This
observation is presented in Theorem \ref{majorr}.

\vspace{.4in}
\section{From $L^2$ to $L^\infty$}
\setcounter{equation}{0} \label{l2l8}

\vspace{.1in} In this section, we show that any Leray-Hopf weak
solution of (\ref{QG}) is actually in $L^\infty$ for $t>0$. More
precisely, we have the following theorem.

\begin{thm}\label{l8bound}
Let $\theta_0\in L^2({\Bbb R}^n)$ and let $\theta$ be a
corresponding Leray-Hopf weak solution of (\ref{QG}). That is,
$\theta$ satisfies
\begin{equation}\label{lh}
\theta \in L^\infty([0,\infty), L^2({\Bbb R}^n))\cap L^2([0,\infty);
\mathring{H}^\alpha({\Bbb R}^n)).
\end{equation}
Then, for any $t>0$,
$$
\sup_{{\Bbb R}^n} |\theta(x,t)| \le
C\,\frac{\|\theta_0\|_{L^2}}{t^\frac{n}{4\alpha}}.
$$
As a special consequence,
$$
\|u(\cdot, t)\|_{BMO({\Bbb R}^n)} \le
C\,\frac{\|\theta_0\|_{L^2}}{t^\frac{n}{4\alpha}}
$$
for any $t>0$.
\end{thm}

This theorem can be proved by following the approach of Caffarelli
and Vasseur \cite{CV}. For the sake of completeness, it is provided
in the appendix.

\vspace{.3in}
\section{The diminishing oscillation result}
\setcounter{equation}{0} \label{lemmas}

This section presents the diminishing oscillation result. We first
recall a theorem of Caffarelli and Silvestre \cite{CS}. It states
that if $L(\theta)$ solves the following initial and boundary value
problem
\begin{equation}\label{ext}
\left\{
\begin{array}{l}
 \nabla\cdot (z^b \nabla L(\theta)) =0, \quad (x,z)\in {\Bbb R}^n
 \times (0,\infty),\\ \\
 L(\theta)(x,0) =\theta(x),\quad x\in {\Bbb R}^n,
\end{array}
\right.
\end{equation}
then
\begin{equation}\label{der} (-\Delta)^\alpha \theta =
\lim_{z\to 0} (-z^b L(\theta)_z)
\end{equation}
where $b= 1-2\alpha$. Furthermore, the boundary-value problem
(\ref{ext}) can be solved through a Poisson formula
$$
L(\theta)(x,z) = P(x, z)\ast \theta \equiv \int_{{\Bbb R}^n} P(x-y,
z) \theta(y)\,dy,
$$
where the Poisson kernel
\begin{equation}\label{ker}
P(x,z)=C_{n,b} \frac{z^{1-b}}{(|x|^2+|z|^2)^\frac{n+1-b}{2}}
=C_{n,\alpha} \frac{z^{2\alpha}}{(|x|^2+|z|^2)^\frac{n+2\alpha}{2}}.
\end{equation}
For notational convenience, we shall write
$$
\theta^*(x,z,t) =L(\theta(\cdot,t))(x,z).
$$
The following notation will be used throughout the rest of the
sections:
$$
f_+ = \max(0, f), \quad B_r \equiv [-r,r]^n \subset {\Bbb R}^n,
\quad Q_r \equiv  B_r \times [0,r] \subset {\Bbb R}^n\times\{t\ge
0\}
$$
and
$$ B^*_r \equiv B_r\times [0,r] \subset {\Bbb R}^n\times {\Bbb R}^+,
\quad Q_r^* \equiv  [-r,r]^n \times [0,r]\times[0,r] \subset {\Bbb
R}^n\times {\Bbb R}^+\times\{t\ge 0\}.
$$

\vspace{.1in}
\begin{thm}\label{dios}
Let $\theta$ be a weak solution to (\ref{QG}) satisfying
$$
\theta \in L^\infty([0,\infty), L^2({\Bbb R}^n))\cap L^2([0,\infty);
\mathring{H}^\alpha({\Bbb R}^n))
$$
with $u$ satisfying (\ref{conu2}) below. Assume
$$
|\theta^*| \le 2 \quad \mbox{in $Q^*_4$}.
$$
Then there exists a $\lambda^*>0$ such that
\begin{equation}\label{lessl}
\sup_{Q^*_1}\theta^* -\inf _{Q^*_1}\theta^*\le 4-\lambda^*.
\end{equation}
\end{thm}

\vspace{.1in} The proof of this theorem relies on three propositions
stated below and will be provided in the appendix. It can be seen
from the proofs of this theorem and related propositions that
$\lambda^*$ may depend on $\|u\|_{L^{\frac{n}{\alpha}}}$ in the
fashion $\lambda^*\,\,\sim\,\exp(-\|u\|^m_{L^{\frac{n}{\alpha}}})$
for some constant $m$.

\vspace{.1in} The first proposition derives a local energy
inequality which bounds the $L^2$-norm of the gradient of $\theta^*$
in terms of the local $L^2$-norms of $\theta$ and $\theta^*$.
\begin{prop}\label{localen}
Let $0<t_1<t_2<\infty$. Let $\theta$ be a solution of (\ref{QG})
satisfying
$$
\theta \in L^\infty([t_1,t_2];L^2({\Bbb R}^n)) \cap L^2([t_1,t_2];
\mathring{H}^{\alpha}({\Bbb R}^n)).
$$
Assume the velocity $u$ satisfies
\begin{equation}\label{conu} u \in
L^\infty([t_1,t_2]; L^{\frac{n}{\alpha}}({\Bbb R}^n)).
\end{equation}
Then, for any cutoff function $\eta$ compactly supported in $B^*_r$
with $r>0$,
$$
\int_{t_1}^{t_2} \int_{B^*_r} z^b\,|\nabla (\eta \theta^*_+)|^2
\,dxdzdt + \int_{B_r} (\eta \theta_+)^2 (t_2,x)\,dx \le \int_{B_r}
(\eta \theta_+)^2 (t_1,x)\,dx
$$
\begin{equation}\label{goat}
+ \,C_1\, \int_{t_1}^{t_2} \int_{B_r} (|\nabla \eta|
\theta_+)^2\,dxdt +\int_{t_1}^{t_2} \int_{B^*_r} z^b(|\nabla \eta|
\theta^*_+)^2 \,dxdzdt,
\end{equation}
where
\begin{equation}\label{const1}
C_1=\|u\|_{L^\infty([t_1,t_2]; \,L^{\frac{n}{\alpha}}({\Bbb R}^n))}.
\end{equation}
If, instead of (\ref{conu}), we assume
\begin{equation}\label{conu2} u \in
L^\infty([t_1,t_2]; C^{1-2\alpha}({\Bbb R}^n)) \quad\mbox{and}\quad
\int_{B_r} u(x,t)\,dx =0,
\end{equation}
then the same local energy inequality (\ref{goat}) holds with $C_1$
in (\ref{const1}) replaced by
\begin{equation}\label{const2}
C_2=\|u\|_{L^\infty([t_1,t_2];\, C^{1-2\alpha}({\Bbb R}^n))}.
\end{equation}
\end{prop}

\vspace{.1in} The following proposition establishes the diminishing
oscillation for $\theta^*$ under the condition that the local
$L^2$-norms of $\theta$ and $\theta^*$ are small.
\begin{prop} \label{tech1}
Let $\theta$ be a solution of the supercritical QG equation
(\ref{QG}) satisfying
$$
\theta\in L^\infty([0,\infty);L^2)\cap
L^2([0,\infty);\mathring{H}^\alpha).
$$
Assume that $u$ satisfies the condition in (\ref{conu2}) and
$$
\theta^* \le 2\quad \mbox{in}\,\,  B^*_4\times[-4,0].
$$
There exist $\epsilon_0>0$ and $\lambda>0$ such that if
\begin{equation}\label{as}
\int_{-4}^0\int_{B^*_4} (\theta^*_+)^2 \,z^b\, dxdzds +
\int_{-4}^0\int_{B_4} (\theta_+)^2 \, dxds \le \epsilon_0,
\end{equation}
then
\begin{equation}\label{deed}
\theta_+ \le 2-\lambda \quad \mbox{on} \,\, B_1\times[-1,0].
\end{equation}
\end{prop}

\vspace{.1in} The proof is obtained by following Caffarelli and
Vasseur and will be presented in the appendix. The following
proposition supplies a condition that guarantees the smallness of
the local $L^2$-norms of $\theta$ and $\theta^*$.
\begin{prop} \label{tech2}
Let $\theta$ be a Leray-Hopf weak solution to the supercritical
equation (\ref{QG}) with $u$ satisfying (\ref{conu2}). Assume that
$$
\theta^* \le 2 \qquad \,\mbox{in $Q_4^*$}
$$
and
$$
|\{(x,z,t)\in Q^*_4 : \,\, \theta^* \le 0\}|_w \ge
\frac{|Q^*_4|_w}{2},
$$
where $|Q^*_4|_w$ denotes the weighted measure of $Q^*_4$ with
respect to $z^b\,dxdzdt$. For every $\epsilon_1>0$, there exists a
constant $\delta_1>0$ such that if
$$
|\{(x,z,t)\in Q^*_4 : \,\,0< \theta^*(x,z,t) <1\}|_w \le \delta_1,
$$
then
$$
\int_{Q_1} \theta^2_+\,dxdt + \int_{Q^*_1} (\theta^*_+)^2
\,z^b\,dx\,dzdt \le C\,\epsilon_1^\alpha.
$$
\end{prop}

The proof of this proposition involves a weighted version of De
Giorgi's isoperimetric inequality. More details will be given in the
appendix. The isoperimetric inequality with no weight was given in
Caffarelli and Vasseur \cite{CV}.
\begin{lemma}\label{degio}
Let $B_r=[-r,r]^n\subset {\Bbb R}^n$ and $B^*_r=B_r\times [0,r]$.
Let $b\in [0,1)$ and let $p>(1+b)/(1-b)$. Let $f$ be a function
defined in $B_r^*$ such that
$$
K\equiv \int_{B_r}\int_0^r z^b |\nabla f|^2\,dz \,dx <\infty.
$$
Let
\begin{equation}\label{abcd}
\begin{array}{l}
\displaystyle {\cal A}\equiv \{(x,z)\in B_r^*: \,f(x,z) \le 0\}, \\
\displaystyle {\cal B}\equiv \{(x,z)\in B_r^*: \,f(x,z) \ge 1\},\\
\displaystyle {\cal C}\equiv \{(x,z)\in B_r^*: 0<\,f(x,z) < 1\}
\end{array}
\end{equation}
and let $|{\cal A}|_w$, $|{\cal B}|_w$ and $|{\cal C}|_w$ be the
weighted measure of ${\cal A}$, ${\cal B}$ and ${\cal C}$ with
respect to $z^b dxdz$, respectively. Then
$$
|{\cal A}|_w\,|{\cal B}|_w \le C
r^{1+\frac{1}{2}(n+1-\frac{p+1}{p-1} b)(1-\frac1p)}\, (|{\cal
C}|_w)^{\frac{1}{2p}} \, K^{\frac12}
$$
where $C$ is a constant independent of $r$.
\end{lemma}
{\it Proof}.\quad We scale the $z-$variable by
$$
\tilde{z} =\frac{1}{b+1} z^{b+1}\quad\mbox{or}\quad
z=((b+1)\tilde{z})^{\frac{1}{b+1}}.
$$
When $(x,z)\in B_r\times [0,r]$, $(x,\tilde{z})\in B_r\times [0,
\tilde{r}]$ with $ \tilde{r} = \frac{r}{1+b}$. For notational
convenience, we write $E_r = B_r \times [0,\tilde{r}]$. Define
$$
g(x,\tilde{z}) =f(x,z)\quad \mbox{for} \quad (x,\tilde{z})\in
B_r\times [0, \tilde{r}].
$$
Let
$$
\tilde{{\cal A}} \equiv \{(x,\tilde{z})\in B_r\times [0, \tilde{r}]:
\,g(x,\tilde{z})\le 0\}
$$
and $\tilde{{\cal B}}$ and $\tilde{{\cal C}}$ be similarly defined.
Therefore,
\begin{eqnarray}
|{\cal A}|_w\,|{\cal B}|_w &\equiv& \int_{{\cal A}} \int_{{\cal
B}}\,z_1^b dx_1dz_1\,
z_2^b dx_2dz_2 \nonumber\\
&\le& \int_{{\cal A}} \int_{{\cal B}} (f(x_1,z_1)-f(x_2,z_2)\,z_1^b
dx_1dz_1\, z_2^b dx_2dz_2 \nonumber\\
&=& \int_{\tilde{{\cal A}}} \int_{\tilde{{\cal B}}}
(g(x_1,\tilde{z_1})-g(x_2,\tilde{z_2}))\,dx_1 d\tilde{z_1}\,dx_2
d\tilde{z_2} \nonumber\\
&=& \int_{\tilde{{\cal A}}} \int_{\tilde{{\cal B}}}
(g(\tilde{y}_1)-g(\tilde{y}_2))\,d\tilde{y}_1
d\tilde{y}_2\label{int1},
\end{eqnarray}
where $\tilde{y}_1=(x_1,\tilde{z_1})$ and
$\tilde{y}_2=(x_2,\tilde{z_2})$. This integral now involves no
weight and can be handled similarly as in Caffarelli and Vasseur
\cite{CV}.
\begin{eqnarray}
|{\cal A}|_w\,|{\cal B}|_w &\le& C\, \int_{E_r}\int_{E_r}
\frac{|\nabla g(\tilde{y}_1 +\tilde{y}_2)|}{|\tilde{y}_2|^{n-1}}
\,\chi_{\{\tilde{y}_1 +\tilde{y}_2)\in
\tilde{{\cal C}}\}}\, d \tilde{y}_1\,d \tilde{y}_2 \label{int2}\\
&=& C\, \int_{E_r} \int_{E_r+\{ \tilde{y}_2\}} |\nabla g(\tilde{y})|
\,\chi_{\{\tilde{y}\}} \,d \tilde{y}\, \frac{1}{|\tilde{y}_2|^{n-1}}
d\tilde{y}_2 \nonumber \\
&=& C\, r\, \int_{E_r} |\nabla g(\tilde{y})|
\,\chi_{\{\tilde{y}\in\, \tilde{{\cal C}}\}} \,d \tilde{y},
\nonumber
\end{eqnarray}
where $\chi$ denotes the characteristic function. By the definition
of $g$,
$$
\nabla g(x,\tilde{z}) = (\nabla_x g, \partial_{\tilde{z}} g) =
(\nabla_x f, \partial_z f \frac{\partial z}{\partial \tilde{z}})
=(\nabla_x f, \partial_z f\, z^{-b}).
$$
By substituting back to the $z-$variable and letting $y=(x,z)$, we
have
\begin{eqnarray}
|{\cal A}|_w\,|{\cal B}|_w &\le&C\, r\, \int_{B_r}\int_0^r
\chi_{\{y\in {\cal C}\}}\,
\sqrt{|\nabla_x f|^2 + (\partial_z f)^2 z^{-2b}}\, z^b \,dz\,dx  \nonumber \\
&\le& C\,r\,\Big(\int_{B_r^*} (|\nabla_x f|^2 z^{2b} + (\partial_z
f)^2)\, z^b \,dzdx\Big)^{1/2}\,\Big(\int_{B_r^*}\chi_{\{y\in \,{\cal
C}\}} z^{-b}\,dzdx\Big)^{1/2}.\nonumber
\end{eqnarray}
By H\"{o}lder's inequality,
\begin{eqnarray}
\int_{B_r^*}\chi_{\{y\in \,{\cal C}\}} z^{-b}\,dzdx &\le&
\Big(\int_{{\cal C}} z^b \,dzdx \Big)^{1/p} \,
\Big(\int_{B_r}\int_0^r z^{-\frac{p+1}{p-1} b}\,dzdx\Big)^{1-1/p}\nonumber \\
&=& |{\cal C}|^{1/p}_w \,\, r^{(n+1-\frac{p+1}{p-1}
b)(1-\frac1p)}.\nonumber
\end{eqnarray}
Therefore,
$$
|{\cal A}|_w\,|{\cal B}|_w \le
C\,r^{1+\frac{1}{2}(n+1-\frac{p+1}{p-1} b)(1-\frac1p) }\,\,|{\cal
C}|^{\frac{1}{2p}}_w\, K^{\frac{1}{2}}.
$$
This completes the proof of this lemma.

\vspace{.3in}
\section{H\"{o}lder continuity under the condition $u\in C^{1-2\alpha}$}

\setcounter{equation}{0} \label{l8toholder}

\vspace{.1in}This section proves the following theorem.
\begin{thm}\label{majorr}
Let $\theta$ be a solution of (\ref{QG}) satisfying
$$
\theta\in L^\infty([0,\infty),L^2({\Bbb R}^n)) \cap
L^2([0,\infty);\mathring{H}^{\alpha}({\Bbb R}^n)).
$$
Let $t_0>0$. Assume that
$$
\theta \in L^\infty({\Bbb R}^n \times [t_0,\infty))
$$
and
$$
u\in L^\infty([t_0,\infty); C^{1-2\alpha}({\Bbb R}^n)).
$$
Then $\theta$ is in $C^\delta({\Bbb R}^n\times [t_0,\infty))$ for
some $\delta>0$.
\end{thm}
{\it Proof}.\quad Fix $x\in {\Bbb R}^n$ and $t\in [t_0,\infty)$. We
show $\theta$ is $C^\delta$ at $(x,t)$.  Define
$$
F_0(y,s) =\theta(x+y+x_0(s),\, t+s),
$$
where $x_0(s)$ is the solution to
\begin{eqnarray}
&& \mathring{x}_0(s) =\frac{1}{|B_4|} \int_{x_0(s)+B_4} u(x+y,
t+s)\,dy, \nonumber\\
&& x_0(0)=0. \nonumber
\end{eqnarray}
Note that $x_0(s)$ is uniquely defined from the classical
Cauchy-Lipschitz theorem. Since $\theta$ is bounded in ${\Bbb R}^n
\times [t_0,\infty)$, we can define
\begin{eqnarray}
&&\overline{\theta}^*_0=\frac{4}{\sup_{Q_4^*}F_0^*-\inf_{Q_4^*}F_0^*}
\Big(F_0^*-\frac{\sup_{Q_4^*}F_0^* +\inf_{Q_4^*}F_0^*}{2}\Big),\nonumber\\
&&u_0(y,s)=u(x+y+x_0(s),\, t+s)-\mathring{x}_0(s),\nonumber
\end{eqnarray}
where $F^*_0(y,z,s) =L(F_0(\cdot,s))(y,z)$. Trivially,
$|\overline{\theta}^*_0| \le 2$  and thus $|\overline{\theta}_0| \le
2$. To verify that $(\overline{\theta}_0, u_0)$ solves the
supercritical QG equation (\ref{QG}), it suffices to show that
$(F_0,u_0)$ solves (\ref{QG}). In fact,
\begin{eqnarray}
\partial_s F_0 + u_0\cdot \nabla_y F_0 &=&
\mathring{x}_0(s)\cdot\nabla_x \theta + \partial_t \theta +
(u-\mathring{x}_0(s))\cdot\nabla_x\theta \nonumber\\
&=&\partial_t \theta + u \cdot\nabla\theta_x =-\Lambda^{2\alpha}_x
\theta =-\Lambda^{2\alpha}_y F_0. \nonumber
\end{eqnarray}
In addition, for any $s\ge 0$,
$$\|u_0(\cdot, s)\|_{C^{1-2\alpha}} =\|u(\cdot, t+s)\|_{C^{1-2\alpha}}
\quad\mbox{and}\quad  \int_{B_4} u_0(y,s)dy =0.
$$

Let $\mu>0$ and set for every integer $k>0$
\begin{eqnarray}
&& F_k(y,s)=\mu^{2\alpha-1}\,F_{k-1}\left(\mu
y+\mu^{2\alpha}\,x_k(s),
\mu^{2\alpha} s\right),\nonumber\\
&&
\overline{\theta}^*_k=\frac{4}{\sup_{Q_4^*}F_k^*-\inf_{Q_4^*}F_k^*}
\Big(F_k^*-\frac{\sup_{Q_4^*}F_k^* +\inf_{Q_4^*}F_k^*}{2}\Big),\nonumber\\
&&\mathring{x}_k(s) =\frac{1}{|B_4|} \int_{B_4+
\mu^{2\alpha-1}x_k(s)} u_{k-1}\Big(\mu y,
\mu^{2\alpha} s\Big)\,dy, \nonumber\\
&&x_k(0)=0, \nonumber \\
&& u_k(y,s)=\mu^{2\alpha-1}\,u_{k-1}\Big(\mu
y+\mu^{2\alpha}\,x_k(s),\, \mu^{2\alpha}
s\Big)-\mu^{2\alpha-1}\,\mathring{x}_k(s).\nonumber
\end{eqnarray}
By the construction, $|\overline{\theta}_k|\le 2$ and
\begin{eqnarray}
\|u_k(\cdot,s)\|_{C^{1-2\alpha}}&=&\mu^{2\alpha-1}\,\|u_{k-1}(\mu
\,\cdot + \mu^{2\alpha},
\mu^{2\alpha}s)\|_{C^{1-2\alpha}}\nonumber \\
&\le & \|u_{k-1}(\cdot,
\mu^{2\alpha}s)\|_{C^{1-2\alpha}}\nonumber \\
&\le& \|u_0(\cdot,\mu^{2\alpha k}s)\|_{C^{1-2\alpha}} \nonumber\\
&=& \|u(\cdot, t+\mu^{2\alpha k}s)\|_{C^{1-2\alpha}}.\nonumber
\end{eqnarray}
Furthermore,
$$
\int_{B_4} u_k(y,s) dy =0.
$$
We show inductively that $(\overline{\theta}_k, u_k)$ solves
(\ref{QG}). Assume that $(\overline{\theta}_{k-1}, u_{k-1})$ solves
(\ref{QG}), we show that $(\overline{\theta}_{k}, u_{k})$ solves
(\ref{QG}). It suffices to show that $(F_k, u_k)$ solves (\ref{QG}).
By construction, we have
\begin{eqnarray}
\partial_s F_k + u_k\cdot \nabla_y F_k &=& \mu^{4\alpha-1}\,
\mathring{x}_k(s)\cdot\nabla F_{k-1} + \mu^{4\alpha-1}\,\partial_s
F_{k-1} \nonumber\\ &&\,+ \mu^{4\alpha-1}\,
(u_{k-1}-\mathring{x}_k(s))\cdot\nabla F_{k-1} \nonumber\\
&=& \mu^{4\alpha-1}\,\left(\partial_s F_{k-1} + u_{k-1} \cdot\nabla
F_{k-1}\right) \nonumber\\ &=&- \mu^{4\alpha-1}\,\Lambda^{2\alpha}
F_{k-1}  \nonumber\\ &=&-\Lambda^{2\alpha}_y F_k. \nonumber
\end{eqnarray}

\vspace{.1in} For every $k$, we apply the diminishing oscillation
result (Theorem \ref{dios}). There exists a $\lambda^*$ such that
$$
\sup_{Q_1^*} \overline{\theta}_k^* -\inf_{Q_1^*}
\overline{\theta}_k^* \le 4-\lambda^*.
$$
$\lambda^*$ is independent of $k$ since $\|u_k\|_{C^{1-2\alpha}}$
obeys a uniform bound in $k$. According to the construction of
$\overline{\theta}_k^*$, we have
$$
\sup_{Q_1^*} \overline{\theta}_k^* -\inf_{Q_1^*}
\overline{\theta}_k^* =
\frac{4}{\sup_{Q_4^*}F_k^*-\inf_{Q_4^*}F_k^*} (\sup_{Q_1^*} F_k^*
-\inf_{Q_1^*} F_k^*).
$$
Therefore,
$$
\sup_{Q_1^*} F_k^* -\inf_{Q_1^*} F_k^* \le
\Big(1-\frac{\lambda^*}{4}\Big)(\sup_{Q_4^*}F_k^*-\inf_{Q_4^*}F_k^*).
$$
By the construction of $F_k$, we have
$$
\sup_{(y,s)\in Q_4^*}F_k^*(y,s)-\inf_{(y,s)\in Q_4^*}F_k^*(y,s)
$$ $$
=\mu^{2\alpha-1} \left(\sup_{(y,s)\in Q_4^*}F^*_{k-1}(\mu y +
\mu^{2\alpha} x_k(s), \mu^{2\alpha} s)-\inf_{(y,s)\in
Q_4^*}F^*_{k-1}(\mu y + \mu^{2\alpha} x_k(s), \mu^{2\alpha}
s)\right).
$$
For notational convenience, we have omitted the $z$-variable. It is
easy to see from the construction of $\mathring{x}_k$ that
\begin{equation}\label{xdot}
|\mathring{x}_k(s)| \le \|u_{k-1}(\cdot,
\mu^{2\alpha}s)\|_{L^\infty} \le \|u_{k-1}(\cdot,
\mu^{2\alpha}s)\|_{C^{1-2\alpha}} \le \|u(\cdot, t+ \mu^{2\alpha
k}s)\|_{C^{1-2\alpha}}.
\end{equation}
For $0\le s \le 1$, we can choose $\mu>0$ sufficiently small such
that
\begin{equation}\label{keysss}
|\mu y + \mu^{2\alpha} x_k(s)| \le 4\mu + C\,\mu^{2\alpha} <1.
\end{equation}
We then have
$$
\sup_{(y,s)\in Q_4^*}F^*_{k-1}(\mu y + \mu^{2\alpha} x_k(s),
\mu^{2\alpha} s)-\inf_{(y,s)\in Q_4^*}F^*_{k-1}(\mu y +
\mu^{2\alpha} x_k(s), \mu^{2\alpha} s)
$$ $$
\le \sup_{(y,s)\in Q_1^*}F^*_{k-1}(y,s) -\inf_{(y,s)\in
Q_1^*}F_{k-1}^*(y,s).
$$
Consequently,
$$
\sup_{Q_1^*} F_k^* -\inf_{Q_1^*} F_k^* \le \mu^{2\alpha-1}
\Big(1-\frac{\lambda^*}{4}\Big)(\sup_{Q_1^*}F^*_{k-1}
-\inf_{Q_1^*}F_{k-1}^*).
$$
By iteration, for any $k>0$,
\begin{equation}\label{base}
\sup_{Q_1^*} F_k^* -\inf_{Q_1^*} F_k^* \le \mu^{(2\alpha-1)k}
\Big(1-\frac{\lambda^*}{4}\Big)^k (\sup_{Q_1^*} F_0^* -\inf_{Q_1^*}
F_0^*).
\end{equation}
By construction,
\begin{eqnarray}
F_0(y,s) &=& \theta(x+y+x_0(s), t+s), \nonumber\\
F_k (y,s) &=& \mu^{(2\alpha-1)k}\,\theta \Big(x+ \mu^k y +
\mu^{2\alpha+k-1} x_k(s) +\mu^{2\alpha+k-2} x_{k-1}(\mu^{2\alpha}s)
\nonumber\\
&& +\cdots + \mu^{2\alpha} x_1(\mu^{2\alpha(k-1)}s)
+x_0(\mu^{2\alpha k} s), t+ \mu^{2\alpha k} s\Big).\nonumber
\end{eqnarray}
To deduce the H\"{o}lder continuity of $\theta$ in $x$, we set
$s=0$. Then (\ref{base}) implies
$$
\sup_{y\in B_1} \mu^{(2\alpha-1)k}\,\theta(x+\mu^k y,t) -\inf_{y\in
B_1}\mu^{(2\alpha-1)k}\,\theta(x+\mu^k y,t) \le C\,
\mu^{(2\alpha-1)k} \Big(1-\frac{\lambda^*}{4}\Big)^k.
$$
or
\begin{equation}\label{goo}
\sup_{y\in B_1} \theta(x+\mu^k y,t) -\inf_{y\in B_1}\theta(x+\mu^k
y,t) \le C\,\Big(1-\frac{\lambda^*}{4}\Big)^k.
\end{equation}
To see the H\"{o}lder continuity from this inequality, we choose
$\delta>0$ such that
$$
1-\frac{\lambda^*}{4} < \mu^\delta.
$$
Then, for any $|y|>0$, we choose $k$ such that
$$
\left(\frac{1-\frac{\lambda^*}{4}}{\mu^\delta}\right)^k \le
|y|^\delta \quad \mbox{or}\quad \Big(1-\frac{\lambda^*}{4}\Big)^k\le
(\mu^k |y|)^\delta.
$$
It then follows from (\ref{goo}) that
$$
\sup_{y\in B_1} \theta(x+\mu^k y,t) -\inf_{y\in B_1}\theta(x+\mu^k
y,t) \le C\, (\mu^k |y|)^\delta.
$$
For general $0\le s\le 1$ and $y\in B_1$, we have, according to
(\ref{xdot}),
\begin{eqnarray}
 r_k &\equiv& \mu^{2\alpha+k-1} x_k(s) +\mu^{2\alpha+k-2}
x_{k-1}(\mu^{2\alpha}s) +\cdots + \mu^{2\alpha}
x_1(\mu^{2\alpha(k-1)}s) +x_0(\mu^{2\alpha k} s) \nonumber \\
&\le& C \mu^{2\alpha+k-1}|s| \,(1+ \mu^{2\alpha-1} +\cdots +
\mu^{(2\alpha-1)k})\nonumber \\
&=& C\, |s|\,\mu^{2\alpha(k+1)-1}
\frac{1-\mu^{(1-2\alpha)(k+1)}}{1-\mu^{1-2\alpha}} \nonumber \\
&\le& C\,|s| \,\mu^{2\alpha(k+1)-1}.\nonumber
\end{eqnarray}
Without loss of generality, we can assume that $\mu^k |y| > |s|
\mu^{2\alpha k}$. Then we can pick up $\delta>0$ satisfying
$$
1-\frac{\lambda^*}{4} < \mu^{2\alpha \delta}
$$
and suitable $k$ such that
$$
\sup_{(y,s)\in B_1\times[0,1]} \theta(x+\mu^k y + r_k,t+\mu^{2\alpha
k}s) -\inf_{(y,s)\in
B_1\times[0,1]}\theta(x+\mu^k\,y+r_k,t+\mu^{2\alpha k}s)
$$
$$ \le C\, (\mu^k |y|)^\delta +C\,
(\mu^{2\alpha k}|s|)^\delta.
$$
That is, $\theta$ is H\"{o}lder continuous at $(x,t)$. This
completes the proof.

\vspace{.4in}\noindent {\bf Acknowledgment}: PC was partially
supported by NSF-DMS 0504213. JW thanks the Department of
Mathematics at the University of Chicago for its support
and hospitality.

\vspace{.6in} \noindent{\large \bf Appendix}

\renewcommand{\thesection}{\Alph{section}}
\renewcommand{\theequation}{\Alph{section}.\arabic{equation}}
\setcounter{section}{1} \setcounter{equation}{0} \setcounter{thm}{0}

\vspace{.15in}

The appendix contains the proofs of several theorems and
propositions presented in the previous sections. These proofs are
obtained by following the ideas of Caffarelli and Vasseur \cite{CV}.
They are attached here for the sake of completeness.

\vspace{.1in}\noindent  {\it Proof of Theorem \ref{l8bound}}.\quad
We first remark that (\ref{lh}) implies that $\theta$ satisfies the
level set energy inequality. That is, for every $\lambda>0$,
$\theta_\lambda=(\theta-\lambda)_+$ satisfies
\begin{equation} \label{level}
\int \theta^2_\lambda(x,t_2) \,dx + 2\int_{t_1}^{t_2} \int
|\Lambda^\alpha\theta_\lambda|^2 dx\,dt \le \int
\theta^2_\lambda(x,t_1) \,dx
\end{equation}
for any $0<t_1<t_2<\infty$. This can be verified by using an
inequality of A. C\'{o}rdoba  and D. C\'{o}rdoba \cite{CC} for
fractional derivatives, namely
$$
f'(\theta) (-\Delta)^\alpha \theta \ge   (-\Delta)^\alpha f(\theta)
$$
for any convex function $f$. Applying this inequality with
$$
f(\theta) = (\theta-\lambda)_+,
$$
we have
$$
\partial_t \theta_\lambda + u\cdot \nabla \theta_\lambda +
\Lambda^{2\alpha} \theta_\lambda \le 0.
$$
Multiplying this equation by $\theta_\lambda$ then leads to
(\ref{level}). Let $k\ge 0$ be an integer and let
$\lambda=C_k=M(1-2^{-k})$ for some $M$ to be determined. It then
follows from (\ref{level}) that
$$
\theta_k = (\theta-C_k)_+.
$$
satisfies
\begin{equation}\label{diffin}
\partial_t \int \theta_k^2(x,t) \,dx + \int |\Lambda^\alpha\theta_k|^2
dx\, \le 0.
\end{equation}
Fix any $t_0>0$. Let $t_k=t_0(1-2^{-k})$. Consider the quantity
$U_k$,
$$
U_k = \sup_{t\ge t_k} \int \theta^2_k(x,t) \, dx +
2\int_{t_k}^\infty\int |\Lambda^\alpha\theta_k|^2 dx\,dt.
$$
Now let $s\in [t_{k-1},t_k]$. We have from (\ref{diffin}) that for
any $s\le t$,
$$
\int \theta^2_k(x,t) \, dx + 2\int_{s}^t\int
|\Lambda^\alpha\theta_k|^2 dx\,dt \le \int \theta_k^2(x,s) dx
$$
which implies that
$$
\sup_{t\ge t_k} \int \theta^2_k(x,t) \, dx \le \int \theta_k^2(x,s)
dx,\quad 2\int_{s}^\infty\int |\Lambda^\alpha\theta_k|^2 dx\,dt \le
\int \theta_k^2(x,s) dx
$$
Since $s\in (t_{k-1},t_k)$, we add up these inequalities to get
$$
U_k \le 2 \int \theta_k^2(x,s) dx.
$$
Taking the mean in $s$ over $[t_{k-1},t_k]$, we get
\begin{equation}\label{ukb}
U_k \le \frac{2^{k+1}}{t_0} \int_{t_{k-1}}^\infty \int
\theta_k^2(x,t) \,dx\,dt
\end{equation}
By Sobolev embedding and Riesz interpolation,
$$
\|\theta_{k-1}\|^2_{L^q([t_{k-1},\infty)\times {\Bbb R}^n)} \le C
\left(\sup_{t\ge t_{k-1}} \int \theta_{k-1}^2(x,t) dx \right)^\sigma
\left(\int_{t_{k-1}}^\infty \int_{{\Bbb R}^n} |\Lambda^\alpha
\theta_{k-1}|^2\,dx\,dt\right)^{1-\sigma},
$$
where
\begin{equation}\label{qde}
\frac{1}{q}=\frac{1-\sigma}{2} =\frac{\sigma}{2} +
\left(\frac{1}{2}-\frac{\alpha}n\right)(1-\sigma), \quad\mbox{or}
\quad \sigma =\frac{2\alpha}{n+2\alpha},\,\,\, q=
2+\frac{4\alpha}{n},
\end{equation}
Therefore,
$$
U_{k-1} \ge C \left(\int_{t_{k-1}}^\infty \int |\theta_{k-1}|^q
\,dx\,dt\right)^{2/q}.
$$

By the definition of $\theta_k$, $\theta_k\ge 0$. When $\theta_k
>0$,
$$
\theta_{k-1} = \theta_k + M2^{-k} \ge M2^{-k}
$$
and thus we have
$$
\chi_{\{(x,t):\,\theta_k>0\}} \le \left(\frac{2^k
\theta_{k-1}}{M}\right)^{q-2},
$$
where $\chi$ denotes the characteristic function. It then follows
from (\ref{ukb}) that
\begin{eqnarray}
U_k &\le & \frac{2^{k+1}}{t_0}\int_{t_{k-1}}^\infty \int
\theta_k^2(x,t) \,\chi_{\{\theta_k>0\}}\, dx\,dt \nonumber \\
&\le& \frac{2^{k+1}}{t_0}\int_{t_{k-1}}^\infty \int
\theta_{k-1}^2(x,t) \,\chi_{\{\theta_k>0\}}\, dx\,dt
\nonumber \\
&\le& \frac{2^{k+1+(q-2)k}}{t_0 M^{q-2}}\,\int_{t_{k-1}}^\infty \int
|\theta_{k-1}|^q \,dx\,dt \nonumber \\
&\le& \frac{2}{t_0 M^{q-2}}\, 2^{(q-1)k}\,
U_{k-1}^{\frac{q}{2}}.\label{basein}
\end{eqnarray}
Since $q>2$, we rewrite (\ref{basein}) as
\begin{equation}\label{vin}
V_k \le V_{k-1}^{\frac{q}2},
\end{equation}
where
$$
V_k = \frac{2^{\gamma k}\,U_k}{t_0^{2/(q-2)}\, M^2 \, 2^{(-\gamma
q-2)/(q-2)}}\quad\mbox{with}\quad \gamma = \frac{2(q-1)}{q-2}>0.
$$
Since $U_0\le \|u_0\|^2_{L^2}<\infty$, we can choose sufficiently
large $M$ such that $V_0<1$ and (\ref{vin}) then implies $V_k\to 0$
as $k\to \infty$. Consequently, we conclude that for each fixed
$t_0>0$ and $M$ sufficiently large, $U_k\to 0$ as $k\to \infty$.
That is, $\theta\le M$. Applying this process to $-\theta$ yields a
lower bound.

\vspace{.1in} The scaling invariance
$$
\theta_\rho(x,t) = \rho^{2\alpha-1}\, \theta(\rho x, \rho^{2\alpha}
t)
$$
of (\ref{QG}) allows us to deduce the following explicit bound
$$
\|\theta(\cdot,t) \|_{L^\infty} \le
C\,\frac{\|u_0\|_{L^2}}{t^{\frac{n}{4\alpha}}}.
$$
This concludes the proof of Theorem \ref{l8bound}.

\vspace{.2in} \noindent{\it Proof of Proposition
\ref{localen}}.\quad Multiplying the first equation in (\ref{ext})
by $\eta^2 \theta^*_+$ and integrating over ${\Bbb R}^n\times
(0,\infty)$ leads to
\begin{eqnarray}
0 &=&\int_0^\infty \int_{{\Bbb R}^n} \eta^2 \theta^*_+ \nabla
\cdot(z^b \nabla \theta^*) dxdz \nonumber \\
&=& \int_0^\infty \int_{{\Bbb R}^n} (\nabla\cdot (\eta^2\theta^*_+\,
z^b \nabla\theta^*)-\nabla(\eta^2\theta^*_+)\cdot z^b
\nabla\theta^*)\,dxdz.\nonumber
\end{eqnarray}
Since $\eta$ has compact support on $B_r^*$ and
$$
\lim_{z\to 0} (-z^b \partial_z \theta^*) = (-\Delta)^\alpha \theta
\equiv \Lambda^{2\alpha} \theta,
$$
we have
\begin{eqnarray}
0 &=&  \int_{{\Bbb R}^n} \eta^2\theta_+ \Lambda^{2\alpha}\theta
dx-\int_0^\infty \int_{{\Bbb R}^n} z^b\,(2\eta \nabla \eta
\theta^*_+\cdot \nabla \theta^* + \eta^2\nabla \theta^*_+\cdot
\nabla \theta^*)\,dxdz \nonumber \\
&=& \int_{{\Bbb R}^n} \eta^2\theta_+ \Lambda^{2\alpha}\theta
dx-\int_0^\infty \int_{{\Bbb R}^n}
z^b\,|\nabla(\eta\theta^*_+)|^2\,dx \,dz \,
\nonumber \\
&& \quad + \int_0^\infty \int_{{\Bbb R}^n} z^b\,|\nabla\eta|^2
(\theta^*_+)^2\,dxdz.\nonumber
\end{eqnarray}
Multiplying both sides of the QG equation (\ref{QG}) by
$\eta^2\,\theta_+$, we get
$$
-\int_{{\Bbb R}^n} \eta^2 \theta_+ \Lambda^{2\alpha}
\theta\,dx=\partial_t \int_{{\Bbb R}^n} \eta^2
\frac{\theta_+^2}{2}\,dx - \int_{{\Bbb R}^n} \nabla(\eta^2) \cdot u
\frac{\theta_+^2}{2}\,dx.
$$
Combining these two equations, we get
$$
\int_0^\infty \int_{{\Bbb R}^n} z^b\,|\nabla(\eta\theta^*_+)|^2\,dx
\,dz \,+ \partial_t \int_{{\Bbb R}^n} \eta^2
\frac{\theta_+^2}{2}\,dx
$$ $$
\qquad =\int_0^\infty \int_{{\Bbb R}^n} z^b\,|\nabla\eta|^2
(\theta^*_+)^2\,dxdz +\int_{{\Bbb R}^n} \nabla(\eta^2) \cdot u
\frac{\theta_+^2}{2}\,dx.
$$
Integrating with respect to $t$ over $[t_1,t_2]$, we get
$$
\int_{t_1}^{t_2}\int_0^\infty \int_{{\Bbb R}^n}
z^b\,|\nabla(\eta\theta^*_+)|^2\,dx \,dz \,dt +\int_{{\Bbb R}^n}
\eta^2 \frac{\theta_+^2}{2}(t_2,x)\,dx
$$
$$
\qquad =\int_{{\Bbb R}^n} \eta^2 \frac{\theta_+^2}{2}(t_1,x)\,dx
+\int_{t_1}^{t_2}\int_0^\infty \int_{{\Bbb R}^n} z^b\,|\nabla\eta|^2
(\theta^*_+)^2\,dxdz
$$
\begin{equation}\label{ter}
+\left|\int_{t_1}^{t_2}\int_{{\Bbb R}^n}\eta\nabla \eta \cdot
u\,\theta_+^2\,dx\,dt\right|.
\end{equation}
We now bound the last term. By the inequalities of H\"{o}lder and
Young,
\begin{equation}\label{fbn}
\left|\int_{{\Bbb R}^n}\eta\nabla \eta \cdot
u\,\theta_+^2\,dx\right| \le \|\eta \theta_+\|_{L^q} \||\nabla\eta|u
\theta_+\|_{L^{q'}} \le \epsilon\,\|\eta \theta_+\|_{L^q}^2 +
\frac{1}{\epsilon} \||\nabla\eta|u \theta_+\|_{L^{q'}}^2,
\end{equation}
where $\epsilon>0$ is small, and $q$ and $q'$ satisfies
$$
\frac{1}{q} = \frac{1}{2} - \frac{\alpha}{n}, \qquad \frac{1}{q} +
\frac{1}{q'} =1.
$$
By the Gagliardo-Nirenberg inequality,
$$
\|\eta \theta_+\|_{L^q}^2 \le C \|\eta \theta_+\|^2_{H^\alpha} =C
\int_{{\Bbb R}^n} \eta \theta_+ \Lambda^{2\alpha} \eta \theta_+\,dx.
$$
Furthermore, since $L(\eta \theta_+)$ and $\eta \theta^*_+$ have the
same trace $\eta \theta_+$ on the boundary $z=0$, we apply the Trace
Theorem to obtain
\begin{eqnarray}
\int_{{\Bbb R}^n} \eta\theta_+\,\Lambda^{2\alpha} \eta \theta_+\,dx
&=& \int_0^\infty \int_{{\Bbb R}^n} z^b |\nabla
(L(\eta\theta_+))|^2\,dx\,dz \nonumber\\
&\le& \int_0^\infty \int_{{\Bbb R}^n} z^b |\nabla
(\eta\theta^*_+)|^2\,dx\,dz.\label{id2}
\end{eqnarray}
Therefore,
\begin{equation}\label{gn}
\|\eta \theta_+\|_{L^q}^2 \le C\int_0^\infty \int_{{\Bbb R}^n} z^b
|\nabla (\eta\theta^*_+)|^2\,dx\,dz.
\end{equation}
Noticing that $1/q'=1/2+\alpha/n$, the second term in (\ref{fbn})
can be bounded by
$$
\||\nabla\eta|u \theta_+\|_{L^{q'}}^2 \le \|u\|_{L^{n/\alpha}}^2
\||\nabla \eta| \theta_+\|_{L^2}^2.
$$
(\ref{goat}) is thus obtained. If we further know that $u$ satisfies
(\ref{conu2}), then
$$
\|u\|_{L^{n/\alpha}} = \left(\int_{B_4}
\left|u(x,t)-\frac{1}{|B_4|}\int_{B_4}
u(y,t)\,dy\right|^{\frac{n}{\alpha}}\,dx \right)^{\frac{\alpha}{n}}
\le C \|u\|_{C^{1-2\alpha}}.
$$
This completes the proof of Proposition \ref{localen}.

\vspace{.2in}\noindent{\it Proof of Theorem \ref{dios}}. It suffices
to show that if
\begin{equation}\label{halfco}
|\{(x,z,t)\in Q^*_4:\,\,\theta^* \le 0\}|_w \ge \frac{1}{2}
|Q^*_4|_w,
\end{equation}
then there exists a $\lambda^*>0$ such that
\begin{equation}\label{halfin}
\theta^* \le 2-\lambda^*\quad\mbox{in $Q^*_1$}.
\end{equation}
Otherwise, we have
$$
|\{(x,z,t)\in Q^*_4:\,\,-\theta^* \le 0\}|_w \ge \frac{1}{2}
|Q^*_4|_w
$$
which implies
$$
-\theta^* \le 2-\lambda^* \quad\mbox{or}\quad \theta^* \ge
-2+\lambda^*\quad\mbox{in $Q^*_1$}.
$$
Thus, in either case,
$$
\sup_{Q_1^*} \overline{\theta}_k^* -\inf_{Q_1^*}
\overline{\theta}_k^* \le 4-\lambda^*.
$$
We now show (\ref{halfin}) under (\ref{halfco}). Fix $\epsilon_0$ as
in (\ref{as}). Choose $\delta_1$ and $\epsilon_1$ as in Proposition
\ref{tech2} with $C\,\epsilon_1^\alpha= \epsilon_0$. Let $K_+$ be
the integer
\begin{equation}\label{k+}
K_+ = \left[\frac{|Q_4^*|_w}{2\delta_1}\right] +1.
\end{equation}
For $k\le K_+$, define
\begin{eqnarray}
&& \overline{\theta}_0 = \theta, \nonumber \\
&& \overline{\theta}_k = 2(\overline{\theta}_{k-1}-1).\nonumber
\end{eqnarray}
It is easy to see that $\overline{\theta}_k = 2^k(\theta-2) +2$.
Note that for every $k$, $\overline{\theta}_k$ verifies (\ref{QG}),
and
$$
\overline{\theta}_k \le 2, \quad \mbox{in $Q_4$},
$$
$$
|\{(x,z,t)\in Q^*_4:\,\,\overline{\theta}^*_k \le 0\}|_w \ge
\frac{1}{2} |Q^*_4|_w.
$$
Assume that for all $k\le K_+$ , $|\{(x,z,t)\in
Q^*_4:\,0<\overline{\theta}^*_k<1\}|_w \ge \delta_1$. Then, for
every $k$,
$$
|\{(x,z,t): \overline{\theta}^*_k \le 0\}|_w =|\{(x,z,t):
\overline{\theta}_{k-1}^* <1\}|_w \ge |\{(x,z,t):
\overline{\theta}_{k-1}^*\le 0\}|_w +\delta_1
$$
Hence,
$$
|\{(x,z,t): \overline{\theta}^*_{K_+} \le  0\}|_w \ge K_+\,\delta_1
+|\{(x,z,t): \theta \le  0\}|_w \ge |Q^*_4|_w.
$$
That is, $\overline{\theta}^*_{K_+}  \le 0$ almost everywhere, which
means
$$
2^{K_+}(\theta^*-2)+2\le 0 \quad\mbox{or}\quad \theta^* \le
2-2^{-K_++1}.
$$
(\ref{lessl}) is then verified by taking $0<\lambda^*<2^{-K_++1}$.

\vspace{.1in} Otherwise, there exists $0\le k_0 \le K_+$ such that
$$
|\{(x,z,t): 0<\overline{\theta}_{k_0}^*<1\}|_w \le \delta_1.
$$
Applying Propositions \ref{tech1} and \ref{tech2}, we get
$\overline{\theta}_{k_0+1} \le 2-\lambda$ which means
$$
\theta \le 2-2^{-(k_0+1)}\lambda \le 2- 2^{-K_+}\lambda\quad
\mbox{in}\,\,Q_2.
$$
Consider the function $f_3$ satisfying
\begin{eqnarray}
&& \nabla\cdot(z^b \nabla  f_3) =0\quad \mbox{in $B_2^*$} \nonumber \\
&& f_3=2 \quad\mbox{on the sides of cube except for $z=0$}
\nonumber \\
&& f_3= 2-2^{-K_+} \inf{(\lambda,1)}\quad \mbox{on $z=0$}. \nonumber
\end{eqnarray}
By the maximum principle, $f_3<2-\lambda^*$ in $B^*_1$ and
$$
\theta^*(x,z,t) \le f_3(x,z,t) < 2-\lambda^* \quad\mbox{in $Q_1^*$}.
$$
This completes the proof of Theorem \ref{dios}.

\vspace{.2in} \noindent{\it Proof of Proposition \ref{tech1}.} We
start with the definition of two barrier functions $f_1$ and $f_2$.
Here $f_1$ satisfies
\begin{equation}\label{bar1}
\left\{
\begin{array}{l}
\nabla\cdot(z^b\nabla f_1) =0, \qquad \mbox{in $B^*_4$},\\
f_1 =2 \qquad \mbox{on the sides of $B_4^*$ except for $z=0$},\\
f_1=0 \qquad\mbox{for $z=0$}.
\end{array}
\right.
\end{equation}
By the maximum principle, for some $\lambda>0$,
$$
f_1(x,z) \le 2-4\lambda\quad\mbox{on $B_2^*$}.
$$
The function $f_2$ satisfies
\begin{equation}\label{bar2}
\left\{
\begin{array}{l}
\nabla\cdot(z^b\nabla f_2)=0 \qquad \mbox{in $[0,\infty)\times
[0,1]$},\\
f_2(0,z) =2 \qquad 0\le z\le 1, \\  f_2(x,0) = f_2(x,1)=0\qquad
0<x<\infty.
\end{array}
\right.
\end{equation}
By separating variables, we can explicitly solve (\ref{bar2}) and
find that
$$
|f_2(x,z)| \le \overline{C}\, e^{-\beta_0 x}
$$
for some constants $\overline{C}>0$ and $\beta_0>0$.

\vspace{.1in} It can be verified that there exist $0<\delta\le 1$
and $M>1$ such that for every $k>0$,
$$
n  \overline{C}\, e^{-\frac{\beta_0}{(2\delta)^k}} \le \lambda
2^{-k-2},\qquad \frac{\|P(\cdot,1)\|_{L^2}}{M^k
\delta^{2\alpha(k+1)}} \le \lambda 2^{-k-2},
$$ $$
C_{0,k}\, M^{-(k-3)(1+\frac{1}{n+1-2\alpha})} \le  M^{-k},\quad
k>12n.
$$
where $P(x,z)$ denotes the Poisson kernel defined in (\ref{ker}) and
$C_{0,k}$ is the constant in (\ref{c0}).

\vspace{.1in} (\ref{deed}) is established through an inductive
procedure, which resembles a local version of the proof for Theorem
\ref{l8bound}. Let $k$ be an integer and set
\begin{equation}\label{thk}
 C_k = 2-\lambda(1+2^{-k}), \qquad \theta_k =(\theta-C_k)_+.
\end{equation}
and let $\eta_k=\eta_k(x)$ be a cutoff function such that
\begin{equation}\label{etk}
\chi_{B_{1+2^{-k-1}}} \le \eta_k \le \chi_{B_{1+2^{-k}}}
\quad\mbox{and}\quad  |\nabla\eta_k| < C 2^k,
\end{equation}
where $\chi$ denotes the characteristic function. Set
\begin{equation}\label{akkk}
A_k = 2 \int_{-1-2^{-k}}^0 \int_{0}^{\delta^k}\int_{{\Bbb R}^n}
z^b\,| \nabla (\eta_k \theta^*_k)|^2 \,dxdzdt + \sup_{[-1-2^{-k},0]}
\int_{{\Bbb R}^n} (\eta_k\theta_k)^2\,dx.
\end{equation}
The goal to prove that
\begin{eqnarray}
 && A_k \le M^{-k},\label{akmk} \\
 && \eta_k \theta^*_k \quad\mbox{is supported in $0\le z\le
\delta^k$}.\label{akck}
\end{eqnarray}
(\ref{deed}) then follows as a consequence of (\ref{akmk}).

\vspace{.1in}We first verify (\ref{akmk}) for $0\le k\le 12n$ and
(\ref{akck}) for $k=0$. Let
$$
T_k =-1-2^{-k} \quad\mbox{and}\quad s\in[T_{k-1},T_k).
$$
Applying (\ref{goat}) with $t_1=s$ and $t_2=t$, we obtain
$$
\int_{s}^{t} \int_{B^*_4} z^b\,|\nabla (\eta \theta^*_+)|^2 \,dxdzdt
+ \int_{B_4} (\eta \theta_+)^2 (t,x)\,dx
$$
$$
\le \int_{B_4} (\eta \theta_+)^2 (s,x)\,dx + C_1 \int_{s}^{t}
\int_{B_4} (|\nabla \eta| \theta_+)^2\,dxdt +\int_{s}^{t}
\int_{B^*_4} z^b(|\nabla \eta| \theta^*_+)^2 \,dxdzdt.
$$
Taking $\sup_{t\in [T_k,0]}$ for both sides and letting $s=T_k$ on
the left gives
$$
\int_{T_k}^0 \int_{B^*_4} z^b\,|\nabla (\eta \theta^*_+)|^2 \,dxdzdt
+  \sup_{t\in [T_k,0]}\int_{B_4} (\eta \theta_+)^2 (t,x)\,dx
$$
$$
\le \int_{B_4} (\eta \theta_+)^2 (s,x)\,dx + C_1 \int_{s}^{0}
\int_{B_4} (|\nabla \eta| \theta_+)^2\,dxdt +\int_{s}^{0}
\int_{B^*_4} z^b(|\nabla \eta| \theta^*_+)^2 \,dxdzdt
$$
$$
\le \int_{B_4} (\eta \theta_+)^2 (s,x)\,dx + C_1 \int_{T_{k-1}}^{0}
\int_{B_4} (|\nabla \eta| \theta_+)^2\,dxdt +\int_{T_{k-1}}^{0}
\int_{B^*_4} z^b(|\nabla \eta| \theta^*_+)^2 \,dxdzdt.
$$
Taking the mean of this inequality in $s$ over $[T_{k-1},T_k]$
yields
\begin{eqnarray}
&& \int_{T_k}^0 \int_{B^*_4} z^b\,|\nabla (\eta \theta^*_+)|^2
\,dxdzdt +  \sup_{t\in [T_k,0]}\int_{B_4} (\eta \theta_+)^2
(t,x)\,dx \nonumber \\
&&\quad  \le 2^k \int_{T_{k-1}}^{T_k}\int_{B_4} (\eta \theta_+)^2
(s,x)\,dx\,ds + C_1 \int_{T_{k-1}}^{0} \int_{B_4} (|\nabla \eta|
\theta_+)^2\,dxdt \nonumber \\
&&\qquad +\int_{T_{k-1}}^{0} \int_{B^*_4} z^b(|\nabla \eta|
\theta^*_+)^2. \,dxdzdt.\label{inte}
\end{eqnarray}
Letting $\eta =\eta_k(x)\phi_k(z)$ with $\phi_k$ supported on
$[0,\delta^k]$ and using the assumption (\ref{as}), we then verify
(\ref{akmk}) for $0<k<12n$ if $\epsilon_0$ satisfies
$$
C2^{24n} (1+C_1)\epsilon_0 \le M^{-12n}.
$$
We now show (\ref{akck}) for $k=0$. By the maximum principle,
$$
\theta^* \le (\theta_+ 1_{B_4})\ast P(\cdot,z) + f_1(x,z)
$$
in $B_4^*\times(0,\infty)$. By construction, the function on the
right-hand side satisfies
$$
\nabla \cdot \left(z^b \nabla ((\theta_+ 1_{B_4})\ast P(z) +
f_1(x,z))\right)=0
$$
and has boundary data greater than or equal to the corresponding
ones for $\theta^*$. To obtain an upper bound for $\theta^*$, we
first notice that $f_1(x,z)\le 2-4\lambda$. In addition,
$$
\|(\theta_+ 1_{B_4})\ast P(\cdot,z)\|_{L^\infty(\{x\in B_4, z\ge
1\})} \le C \|P(\cdot, 1)\|_{L^2}\sqrt{\epsilon_0} \le C
\sqrt{\epsilon_0}.
$$
Here we used $\|\theta_+ 1_{B_4}\|_{L^2}\le C \sqrt{\epsilon_0}$,
which can be deduced from (\ref{as}) through a simple argument.
Choose $\epsilon_0$ small enough to get
$$
\theta^* \le 2-2\lambda \qquad \mbox{for $z\ge 1,\,t\le 0$ and $x\in
B_4$}.
$$
Therefore,
$$
\theta_0^* =(\theta^*-(2-2\lambda))_+ \le 0\qquad \mbox{for $z\ge
1,\,t\le 0$ and $x\in B_4$}.
$$
Hence, $\eta_0 \theta^*_0$ is supported in $0\le z\le \delta^0=1$.

\vspace{.1in} Now, assuming that (\ref{akmk}) and (\ref{akck}) are
verified at $k$, we show they are also true at $k+1$. In the
process, we will also show for each $k$,
\begin{equation}\label{conv}
\eta_k \theta_{k+1}^* \le [(\eta_k \theta_k)\ast P(z)]\,\eta_k
\end{equation}
in the set $\bar{B}_k^*=B_{1+2^{-k}}\times[0,\delta^k]$. First we
control $\theta^*_k$ in $\bar{B}_k^*$ by a function $f$ satisfying
$$
\nabla\cdot(z^b \nabla f)) =0
$$
by considering the contributions on the boundaries. No contributions
come from $z=\delta^k$ thanks to the induction property on $k$. The
contribution from $z=0$ can be controlled by $\eta_k\theta_k\ast
P(\cdot, z)$ since it has the same boundary data as $\theta^*_k$ on
$B_{1+2^{-k-1}}$. On each of the other sides, the contribution can
be controlled by
$$
f_2((-x_i+x^+)/\delta^k, z/\delta^k) + f_2((x_i-x^-)/\delta^k
,z/\delta^k),
$$
where $x^+=1+2^{-k}$ and $x^-=-x^+$. Recall that $f_2$ satisfies
$\nabla\cdot(z^b\nabla f_2)=0$ and is no less than $2$ on the sides
$x_i^+$ and $x_i^-$. By the maximum principle,
$$
\theta^*_k \le \sum_{i=1}^n [f_2((x_i-x^+)/\delta^k, z/\delta^k) +\,
f_2((-x_i+x^-)/\delta^k, z/\delta^k)] + (\eta_k\theta_k)\ast
P(\cdot,z).
$$
We know that, for any $x\in B_{1+2^{-k}}$,
$$
\sum_{i=1}^n [f_2((-x_i+x^+)/\delta^k, z/\delta^k)
+f_2((x_i-x^-)/\delta^k, z/\delta^k)] \le
n\overline{C}e^{-\frac{\beta_0}{(2\delta)^k}} \le \lambda 2^{-k-2}.
$$
Therefore,
$$
\theta^*_k \le (\eta_k\theta_k)\ast P(z) + \lambda 2^{-k-2}.
$$
Consequently,
$$
\theta^*_{k+1} \le (\theta^*_k-\lambda 2^{-k-1})_+ \le
((\eta_k\theta_k)\ast P(z)-\lambda 2^{-k-2})_+
$$
Since, for $z=\delta^{k+1}$,
$$
|(\eta_k\theta_k)\ast P(\cdot,z)| \le A_k \|P(\cdot,z)\|_{L^2} \le
\frac{M^{-k}}{\delta^{2\alpha(k+1)}} \|P(\cdot,1)\|_{L^2} \le
\lambda 2^{-k-2},
$$
we obtain
$$
\eta_{k+1}\theta_{k+1}^* \le 0 \quad\mbox{on}\quad z=\delta^{k+1}.
$$

\vspace{.1in} Let $k>12n +1$. Assuming that (\ref{akmk}) is true for
$k-3$, $k-2$ and  $k-1$, we show
\begin{equation}\label{ak}
A_k \le C_{0,k}\,A_{k-3}^{1+\frac{1}{n+1-2\alpha}},
\end{equation}
where
\begin{equation}\label{c0}
C_{0,k} =\frac{C
\,2^{(1+\frac{4\alpha}{n+1-2\alpha})k}}{\lambda^{\frac{2\alpha}{n+1-2\alpha}}}
.
\end{equation}
Since $\eta\theta^*_+$ has the same boundary condition at $z=0$ as
$(\eta\theta_+)^*$,
$$
\int_0^\infty \int_{{\Bbb R}^n} z^b\,|\nabla(\eta \theta_+^*)|^2\,
dxdz \ge \int_0^\infty \int_{{\Bbb R}^n} z^b\, |\nabla
(\eta\theta_+)^*|^2\,dxdz =\int_{{\Bbb R}^n} |\Lambda^\alpha (\eta
\theta_+)|^2\,dx.
$$
Letting $\eta=\eta_k(x)$ and integrating with respect to $t$ over
$[-1-2^{-k},0]$, we obtain
$$
\int_{-1-2^{-k}}^0 \int_0^{\delta^k} \int_{{\Bbb R}^n}
z^b\,|\nabla(\eta_k\theta_+^*)|^2\,dxdzdt \ge
\int_{-1-2^{-k}}^0\int_{{\Bbb R}^n}|\Lambda^\alpha (\eta
\theta_+)|^2\,dxdt.
$$
According to the definition of $A_k$ in (\ref{akkk}),
$$
A_{k-3} \ge \int_{-1-2^{-k+3}}^0\int_{{\Bbb R}^n}|\Lambda^\alpha
(\eta_{k-3} \theta_{k-3})|^2\,dxdt.
$$
By the Gagliardo-Nirenberg inequality
$$
A_{k-3} \ge C\|\eta_{k-3}\theta_{k-3}\|^2_{L^q([-1-2^{-k+3},0]\times
{\Bbb R}^n)},
$$
where $q$ is defined in (\ref{qde}), namely
$$
q= 2+\frac{4\alpha}{n}.
$$
It then follows from (\ref{conv}) that
$$
\|\eta_{k-3} \theta^*_{k-2}\|^2_{L^q} \le \|P(\cdot,1)\|_{L^1}^2
\|\eta_{k-3}\theta_{k-3}\|^2_{L^q}
$$
Therefore,
$$
A_{k-3} \ge C\|\eta_{k-3}\theta^*_{k-2}\|_{L^q}^2 + C
\|\eta_{k-3}\theta_{k-3}\|^2_{L^q}
$$
$$
\qquad \ge C (\|\eta_{k-1}\theta^*_{k-1}\|_{L^q}^2  +
\|\eta_{k-1}\theta_{k-1}\|^2_{L^q})
$$
The second inequality above follows from the simple fact that
$$
\theta_{k-3} \ge \theta_{k-1} \quad\mbox{and}\quad \eta_{k-3}\ge
\eta_{k-1}.
$$
Letting $\eta =\eta_k(x)$ in (\ref{inte}) yields
$$
A_k \le\,C 2^k (C_1+2) \left(\int \eta_{k-1}^2 \theta_k^2\,dx + \int
\eta_{k-1}^2 (\theta^*_k)^2\,dxdz \right).
$$
The same trick as in the proof of Theorem \ref{l8bound} can then be
played here. If $\theta_k>0$, then $\theta_{k-1}\ge 2^{-k}\lambda$
and thus
$$
\chi_{\{\theta_k>0\}} \le \left(\frac{2^k
\theta_{k-1}}{\lambda}\right)^{q-2} \quad \mbox{and}\quad
\chi_{\{\theta^*_k>0\}} \le \left(\frac{2^k
\theta^*_{k-1}}{\lambda}\right)^{q-2}.
$$
Then,
$$
A_k \le \frac{C \,2^{(q-1)k}}{\lambda^{q-2}} A_{k-3}^{\frac{q}{2}} =
\frac{C
\,2^{(1+\frac{4\alpha}{n+1-2\alpha})k}}{\lambda^{\frac{2\alpha}{n+1-2\alpha}}}
\,A_{k-3}^{1+\frac{1}{n+1-2\alpha}}
=C_{0,k}\,A_{k-3}^{1+\frac{1}{n+1-2\alpha}}.
$$
This completes the proof of Proposition \ref{tech1}.

\vspace{.2in}\noindent{\it Proof of Proposition \ref{tech2}}.\quad
It suffices to show
$$
\int_{Q_1} (\theta-1)^2_+\,dxdt + \int_{Q^*_1} (\theta^*-1)^2_+
\,z^b\,dxdzdt \le C \,\epsilon_1^\alpha.
$$

From the fundamental local energy inequality (\ref{goat}), we have
$$
\int_{-4}^0 \int_{B^*_4} |\nabla \theta^*_+|^2\, z^b\,dx dzdt \le C.
$$
Take $\epsilon_1<<1$ and set
$$
K=\frac{4\int_{-4}^0 \int_{B^*_4}  |\nabla \theta^*_+|^2
z^b\,dxdzdt}{\epsilon_1}.
$$
We further write
$$
I_1= \left\{t\in[-4,0]:\, \int_{B^*_4} |\nabla \theta^*_+|^2(t)\,
z^b\, dxdz\le K \right\}.
$$
It follows from the Chebyshev inequality that
\begin{equation}\label{i1s}
\left|[-4,0]\setminus I_1\right| \le \frac{\epsilon_1}{4}.
\end{equation}
For all $t\in I_1$, the De Giorgi inequality in Lemma \ref{degio}
gives
$$
|{\cal A}(t)|_w |{\cal B}(t)|_w \le C\,|{\cal
C}(t)|_w^{\frac{1}{2p}}\,K^{\frac12},
$$
where ${\cal A}$, ${\cal B}$ and ${\cal C}$ are defined in
(\ref{abcd}) with $r=4$. Set
$$
\delta_1=\epsilon_1^{2p(1+\frac{1}{\alpha})+2},\quad I_2=\{
t\in[-4,0]:\, |{\cal C}(t)|_w^{\frac{1}{2p}} \le
\epsilon_1^{1+\frac{1}{\alpha}}\,\}.
$$
Again by the Chebyshev inequality,
\begin{equation}\label{i2s}
\left|[-4,0]\setminus I_2\right| \le \frac{|\{(x,z,t):\,
0<\theta^*<1\}|_w}{\epsilon^{2p(1+\frac{1}{\alpha})}_1} \le
\frac{\delta_1}{\epsilon_1^{2p(1+\frac{1}{\alpha})} } \le
\epsilon_1^2 \le \frac{\epsilon_1}{4}.
\end{equation}
Now, set $I=I_1\cap I_2$. According to (\ref{i1s}) and (\ref{i2s}),
$$
|[-4,0]\setminus I| \le \frac{\epsilon_1}{4} + \frac{\epsilon_1}{4}
=\frac{\epsilon_1}{2}.
$$
In addition, if $t\in I$ satisfying $|{\cal A}(t)|_w \ge
\frac{1}{4}$, then
\begin{equation}\label{isop}
|{\cal B}(t)|_w \le \frac{C\, |{\cal C}(t)|_w^{\frac1{2p}}\,
K^{\frac12}}{|{\cal A}(t)|_w} \le 4
C\epsilon_1^{\frac12+\frac1{\alpha}}.
\end{equation}
Therefore,
$$
\int_{B_4^*}  (\theta^*_+)^2(t) \,z^b\, dxdz  \le 4\int_{{\cal
B}\cup {\cal C}}z^b \, dxdz \le 4 (|{\cal B}|_w + |{\cal C}|_w) \le
16C\,\epsilon_1^{\frac12+\frac1{\alpha}}.
$$
Let
\begin{equation}\label{p1c}
p_1>\frac{2(1+b)}{1-b}=\frac{2-2\alpha}{\alpha} \quad\mbox{and}
\quad \frac1{p_1} +\frac1{q_1}=\frac12.
\end{equation}
Then $1-(\frac12+\frac1{p_1})bq_1>0$ and by H\"{o}lder's inequality,
\begin{eqnarray}
\int_{B_4}  \theta_+^2(t)  dx &\le&  \int_{B_4} (\max_{z}
\theta^*_+(x,z))^2\,dx \nonumber \\
&\le & 2 \int_{B_4} \int_0^4 |\theta^*|\, |\partial_z
\theta^*|\,dz\,dx\nonumber \\
&= & 2\int_{B_4} \int_0^4
z^{\frac{b}{p_1}}|\theta^*|\,\,z^{\frac{b}{2}}|\partial_z
\theta^*|\,z^{-(\frac12 +\frac{1}{p_1})b}\,dz\,dx \nonumber \\
&\le & 2 \int_{B_4} \Big(\int_0^4 z^b
|\theta^*|^{p_1}\,dz\Big)^{\frac{1}{p_1}} \,\Big(\int_0^4 z^b
|\partial_z \theta^*|^2 \,dz \Big)^{\frac12} \,  \Big(\int_0^4
z^{-(\frac12+\frac1{p_1})bq_1}\,dz\Big)^{\frac{1}{q_1}}\,dx \nonumber \\
&\le& C\,\Big(\int_{B^*_4}
z^b|\theta^*|^{p_1}\,dxdz\Big)^{\frac1{p_1}} \, \Big(\int_{B^*_4}
z^b |\nabla
\theta^*|^2\,dxdz\Big)^{\frac12} \,\nonumber \\
&\le & C\, K^{\frac12}\,\Big(\int_{B^*_4}
z^b\,|\theta^*|^{2}\,dxdz\Big)^{\frac1{p_1}} \nonumber \\
&\le & C\,
\epsilon^{{(\frac12+\frac1{\alpha})\frac{1}{p_1}}-\frac12} \equiv
C\,\epsilon_1^\nu.
\end{eqnarray}
where, thanks to (\ref{p1c}),
$$
 \nu= \Big(\frac12+\frac1{\alpha}\Big)\frac{1}{p_1}-\frac12>0.
$$

\vspace{.1in} The next major part proves that $|{\cal A}(t)|_w \ge
\frac14$ for every $t\in I\cap [-1,0]$. Since
$$
|\{(x,z,t): \, \theta^* \le 0\}|_w \ge \frac{|Q^*_4|_w}{2},
$$
there exists a $t_0\le -1$ such that $|{\cal A}(t_0)|_w
\ge\frac{1}{4}$. Thus, for this $t_0$,
$$
\int_{B_4^*} \theta_+ (t_0)^2 \,dx \le C \epsilon^\nu_1.
$$
Using the local energy inequality (\ref{goat}), we have for all
$t\ge t_0$,
$$ \int_{B_4^*} \theta_+^2(t) dx \le \int_{B_4^*} \theta_+^2(t_0) \,dx +C(t-t_0).
$$
For $t-t_0\le \delta^*=\frac{1}{64C}$, we have
$$
\int_{B_4^*} \theta_+^2(t) \,dx \le \frac{1}{64}.
$$
Since $\delta^*$ does not depend on $\epsilon_1$, we can assume that
$\epsilon_1<<\delta^*$. By
$$
\theta^*_+(x,z,t) \le \theta_+(x,t)+ \int_0^z \partial_z\theta^*_+
dz,
$$
we have
\begin{eqnarray}
z^b\, (\theta_+^*)^2(x,z,t) &\le& 2 z^b \theta^2_+(x,t) + 2z^b
\Big(\int_0^z
\partial_z \theta^*_+ dz \Big)^2 \nonumber\\
&\le & 2 z^b \theta^2_+(x,t) + 2 z \int_0^z z^b |\nabla \theta^*|^2
dz.\nonumber
\end{eqnarray}
For $t-t_0\le \delta^*$, $t\in I$ and $z\le \epsilon_1^2$,
\begin{eqnarray}
\int_0^{\epsilon_1^2} \int_{B_4} z^b\, (\theta_+^*)^2\,dxdz &\le&
\frac{2}{b+1}\epsilon_1^{4-4\alpha} \int_{B_4}\theta^2_+(x,t)\,dx
 + \epsilon_1^4 \int_0^{\epsilon_1^2} \int_{B_4} z^b |\nabla
 \theta^*|^2\,dxdz \nonumber \\
&\le& \frac{1}{64} \epsilon_1^2  + C\, \epsilon^3_1 \le \frac{1}{4}
\epsilon_1^2. \nonumber
\end{eqnarray}
By Chebyshev inequality,
$$
|\{(x,z):\, z\le \epsilon_1^2, x\in B_4, \,\theta^*_+(t) \ge 1\}|_w
\le \frac{\epsilon_1^2}{4}.
$$
Since $|{\cal C}(t)|_w \le \epsilon_1^{2p(1+\frac{1}{\alpha})}$,
this gives
$$
|{\cal A}(t)|_w \ge \epsilon_1^2-
\frac{1}{4}\epsilon_1^2-\epsilon_1^{2p(1+\frac{1}{\alpha})}\ge
\frac{1}{2}\epsilon^2.
$$
Combining this bound with (\ref{isop}) leads to
$$
|{\cal B}(t)|_w \le 4C\, \sqrt{\epsilon_1}.
$$
In turn, this bound leads to
$$
|{\cal A}(t)|_w \ge 1-|{\cal B}(t)|_w-|{\cal C}(t)|_w \ge
1-4C\,\sqrt{\epsilon_1}-\epsilon_1^{2p(1+\frac{1}{\alpha})} \ge
\frac{1}{4}.
$$
Hence, for every $t\in [t_0,t_0+\delta^*]\cap I$, we have $|{\cal
A}(t)|_w \ge \frac{1}{4}.$ On $[t_0+\frac{\delta^*}{2}, t_0+
\delta^*]$, there is $t_1\in I$. The reason is that
$|[-4,0]\setminus I| \le \frac{\epsilon}{2}$ and
$|[t_0+\frac{\delta^*}{2},t_0+\delta^*]| =\frac{\delta^*}{2}
>\frac{\epsilon}{2}$.

\vspace{.1in}This process allows us to construct an increasing
sequence $t_n$, $0\ge t_n\ge t_0+\frac{\delta^*}{2}$ such that
$|{\cal A}(t)|_w \ge \frac{1}{4}$ on $[t_n,t_n+\delta^*]\cap I$.
Since $\delta^*$ is independent of $t_n$, we have
$$
|{\cal A}(t)|_w \ge \frac{1}{4} \quad \mbox{for $t\in I\cap
[-1,0]$}.
$$
According to (\ref{isop}), this gives
$$
|{\cal B}(t)|_w \le 4 C\epsilon_1^{\frac12+\frac1{\alpha}}\le
\frac{\epsilon_1}{16}\quad \mbox{for $t\in I\cap [-1,0]$}.
$$
Therefore,
\begin{eqnarray}
|\{(x,z,t): \,\theta^*\ge 1\}|_w &=& |\{(x,z,t): t\in I\cap [-1,0],
\,\,\theta^*\ge 1\}|_w \nonumber \\
&& + |\{(x,z,t): t\in [-1,0]\setminus I, \,\theta^*\ge 1\}|_w
\nonumber \\
&\le& \frac{\epsilon_1}{16} + \frac{\epsilon_1}{2}\le \epsilon_1.
\nonumber
\end{eqnarray}
Since $(\theta^*-1)_+ \le 1$,
\begin{equation}\label{ugoo}
\int_{Q^*_1} z^b\,(\theta^*-1)_+^2  dxdzdt \le \epsilon_1.
\end{equation}
For fixed $x$ and $t$,
$$
\theta(x,t)-\theta(z)^*(x,z,t) =-\int_0^z \partial_z\theta^* dz.
$$
Thus,
\begin{eqnarray}
z^b\,(\theta-1)^2_+ &\le& 2 z^b(\theta^*(z)-1)^2_+ + z^b
\left(\int_0^z|\nabla\theta^*|\,dz\right)^2  \nonumber \\
&\le& 2 z^b(\theta^*(z)-1)^2_+ +z \int_0^z
z^b|\nabla\theta^*|^2\,dz.\nonumber
\end{eqnarray}
Taking the average in $z$ over $[0,\sqrt{\epsilon_1}]$, we get
$$
\epsilon_1^{\frac{b}{2}}\,(\theta-1)^2_+
\le\frac{2}{\sqrt{\epsilon_1}}\int_0^{\sqrt{\epsilon_1}}z^b\,(\theta^*-1)^2_+
dz +
\sqrt{\epsilon_1}\int_0^{\sqrt{\epsilon_1}}z^b|\nabla\theta^*|^2\,dz.
$$
Integrating with respect to $(x,t)\in B_1\times [0,1]$ and invoking
(\ref{ugoo}) lead to
$$
\int_{Q_1} (\theta-1)^2_+ dxds \le C \epsilon_1^\alpha.
$$

\end{document}